\documentclass[12pt]{article}
\usepackage{amsmath}
\newcommand{\be}{\begin{equation}}
\newcommand{\ee}{\end{equation}}
\newcommand{\ba}{\begin{eqnarray}}
\newcommand{\ea}{\end{eqnarray}}
\newcommand{\ban}{\begin{eqnarray*}}
\newcommand{\ean}{\end{eqnarray*}}

\newcommand{\sect}[1]{\section{#1}  \setcounter{equation}{0}}

\begin{document}
\newtheorem{defn}[lem]{Definition}
\newtheorem{theo}[lem]{Theorem}
\newtheorem{cor}[lem]{Corollary}
\newtheorem{prop}[lem]{Proposition}
\newtheorem{rk}[lem]{Remark}
\newtheorem{ex}[lem]{Example}
\newtheorem{note}[lem]{Note}
\newtheorem{conj}[lem]{Conjecture}

\title{Entropy Functionals, Sobolev Inequalities And $\kappa$-Noncollapsing Estimates Along The Ricci Flow}
\author{Rugang Ye \\ {\small Department  of Mathematics} \\
{\small University of California, Santa  Barbara}}
\date{September 10, 2007}
\maketitle

\noindent 1. Introduction \\
2. Hamilton's entropy functional \\
3. Perelman's entropy functional\\
4. The log entropy functional and the log Sobolev constant\\
5. The logarithmic Sobolev and Sobolev inequalities along the Ricci flow\\
6. The $\kappa$-noncollapsing estimates \\
7. The modified Ricci flow \\
8. Furthe Sobolev inequalities\\
9. The Ricci flow with surgeries

\sect{Introduction}
 
Let $M$ be a smooth manifold of dimension $n \ge 2$.  The Ricci flow 
\ba
\frac{\partial g}{\partial t}=-2Ric
\ea 
starting at a given Riemannian metric $g_0$ on $M$ deforms the Riemannian metric 
$g=g(t)$ in the opposite direction of its Ricci curvature tensor $Ric=Ric_{g(t)}$. 
It turns out that the Ricci flow becomes a nonlinear evolution equation of the parabolic type after 
a suitable gauge fixing is performed, see e.g. [CK].  In other words, we can view the Ricci flow as 
a nonlinear heat equation for Riemannian metrics. Indeed, the evolution equations associated with the Ricci flow for various curvature quantities 
such as the scalar curvature, the Ricci curvature tensor and the Riemann curvature 
tensor are all nonlinear heat equations.  For this reason it is natural to look for geometric quantities 
associated with the Ricci flow which can be motivated by quantities in thermodynamics 
and related theories such as statistical mechanics and information theory.  A central concept in thermodynamics 
is entropy. Entropy, from the Greek $\mu\epsilon\tau\alpha\tau\rho o \pi \acute{\eta}$, meaning 
``transformation", is a measure of the unavailability of a system's energy to do the work. 
In terms of statistical mechanics, the entropy (more precisely, the Boltzmann 
entropy) describes the number of the possible microscopic configurations of a given system. The law of entropy, or the second law of thermodynamics, states that spontaneous changes in isolated 
systems occur with an increase in entropy. 

So far two different kinds of entropy functional have been introduced into the theory of the Ricci flow. 
Both are motivated by concepts of entropy in thermodynamics-statistical mechanics-information theory.
One is R.~Hamilton's entropy, the other is G.~Perelman's entropy.  A major difference between these 
two concepts is this.  In Hamilton's entropy, the scalar curvature $R$ of the metric is viewed as 
the leading quantity of the system and plays the role of a probability density, while in Perelman's entropy 
the leading quantity describing the system is the metric itself.  We'll explain below the motivation of 
Hamilton's entropy in terms of Shannon entropy.  In [P1], Perelman provided 
an interpretation of his entropy in terms of concepts in statistical mechanics including
partition function, energy and entropy. We'll briefly explain this interpretation, and provide a different, very natural motivation for Perelman's entropy in terms       
of the logarithmic Sobolev inequality on the euclidean space ${\bf R}^n$. 

Hamilton established the monotonicity of his entropy along the volume-normalized Ricci flow on the 2-sphere 
$S^2$ [H1]. He and B.~Chow used this monotonicity and 
its generalization to study the Ricci flow on  
$S^2$ [H1] [Ch].  In [P1], Perelman established the monotonicity of his entropy
along the Ricci flow in all dimensions.
As an important application of this entropy monotonicity Perelman derived in [P1]  
the $\kappa$-noncollpasing property of the Ricci flow relative to upper bounds for $|Rm|$, the norm of 
the Riemann curvature tensor, under a finite upper bound for the 
time.    Later, Perelman and the present author 
improved this result independently by replacing the upper bounds for $|Rm|$ with 
upper bounds for the scalar curvature, see  [KL] and [Y1]. Note that the $\kappa$-noncollasping 
property is a crucial ingredient in Perelman's work on the Ricci flow and the geometrization conjecture and 
the Poincar\'{e} conjecture.

  Recently, Perelman's entropy theory for the Ricci flow was further developed in [Y3], [Y4],[Y5], [Y6],
  [Y7] and [Y8].
In [Y3], the logarithmic Sobolev inequality along the Ricci flow was established. 
In general, this inequality depends on a finite upper bound for the time. 
But no such bound is required if the first eigenvalue $\lambda_0(g_0)$ of the operator 
$-\Delta+\frac{R}{4}$ for the initial metric is positive.  
$W^{1,2}$ Sobolev inequalities and $\kappa$-noncollapsing estimates along the Ricci flow 
were then derived in [Y3] as consequences of the logarithmic Sobolev inequality.  
In [Y4] and [Y5], the results of [Y3] were extended to the dimension $n=2$ and the 
case $\lambda_0(g_0)=0$.  In [Y6], the log entropy functional was introduced and its monotonicity along the Ricci flow was 
established, based on Perelman's entropy monotonicity.   As a consequence, it was shown that the logarithmic Sobolev inequality improves along the 
Ricci flow.  
In [Y7], $W^{1,p}$ and $W^{2,p}$ Sobolev inequalities along the Ricci flow were obtained for 
$p \not =2$. Several methods 
were used in [Y7], including some tools from harmonic analysis and potential theory such as  
Bessel potentials and Riesz transforms. 

Finally, Sobolev inequalities and $\kappa$-noncollapsing 
estimates were established for the Ricci flow with surgeries  in [Y8]. The key construction in Perelman's work on the Ricci flow and the geometrization conjecture and the Poincar\'{e} conjecture  is 
the Ricci flow with surgeries [P2], which extends Hamilton's earlier work on surgeries of the Ricci flow 
in a substantial way. The results in [Y8] can be used to replace 
the rather complicated arguments in [P2] for preserving the $\kappa$-noncollapsing property after 
surgeries.   Moreover, the $\kappa$-noncollapsing property is 
established in [Y8] independent of the other properties, making the choice of the surgery 
parameters much simpler. It also becomes easier to establish e.g. the canonical neighborhood 
property. 

We'll review below these recent results, providing helpful clues whenever appropriate. We hope that 
our accounts  can provide an integrated picture of the whole theory, offer some new perspectives, 
and explain the main ideas 
without getting into the technical details of the proofs.

\sect{Hamilton's entropy functional}

Let $M$ be a closed manifold of dimension $n=2$ diffeormorphic to the 2-sphere $S^2$. R.~Hamilton introduced the following entropy functional [H1]
\ba
S_0(g)=\int_M R \ln R dvol
\ea
for metrics $g$ of positive scalar curvature. It can be extended straightforwardly to metrics of nonnegative scalar curvature
because $x \ln x \rightarrow 0$ as $x \rightarrow 0^+$.  This formula of entropy can be motivated by 
Shannon entropy in information theory, which in turn is motivated by Gibbs entropy in statistical mechanics. 

Consider a discret thermodynamic system $X$. The Boltzmann  entropy of $X$ is defined to be 
\ba
S_{B}=k_B \ln W,
\ea
where $k_B$ is the Boltzmann constant and $W$ denotes th number of microstates in $X$. The Gibbs 
entropy of $X$ is defined to be 
\ba
S_G=-k_B \sum_i p_i \ln p_i,
\ea
where $p_i$ denotes the probability of the $i$-th microstate. 
If the mircrostates 
are equiprobable, then $S_B$ is reduced to $S_G$. In general, $S_G$ can be derived from $S_B$ 
as a limit when the number of microstates is very large, see e.g. [CC].    Next let $\mathcal M$ be a 
discret message space with 
a probability measure. 
Then the Shannon entropy of $\mathcal{M}$ is defined to be 
\ba
S=-\sum_i p_i \ln p_i,
\ea
where $p_i$ is the probability of the message $m_i$ taken from the message space $\mathcal M$. 
There is an integral version of Shannon entropy for a continuous message space $\mathcal{M}$, called the differential entropy, or 
differential Shannon entropy, which is defined as follows 
\ba
S_D=-\int_{\mathcal{M}} f \ln f d\mu,
\ea
where $f$ denotes the probability density, and $\mu$ the background measure. (A word of caution: In passing from 
$S$ to $S_D$ in terms of a limit of Riemann sums, one has to throw out an infinite quantity.) 
It is the differential Shannon entropy which can be used to motivate Hamilton's entropy.

Consider a positive solution $u$ of the heat equation 
\ba \label{heat}
\frac{\partial u}{\partial t}=\Delta u
\ea
on $M$, with $\Delta=\Delta_g$ for a given Riemannian metric $g$ on $M$. In the framework of thermodynamics 
$u$ represents the temperature of a physical system with the backgorund given by the Riemannian manifold 
$(M ,g)$.  But we view $M$ as a message space with the evolving probability density function 
$u$ with respect to the background measure given by the volume measure of $g$. Since 
the average of $u$ is preserved by the heat flow \label{heat}, we can normalize $u$  
such that $udvol$ is a a probability measure, i.e. $\int_M u dvol=1$.     The differential Shannon entropy  of $u$ 
 is then
\ba
S_D(u)=-\int_M u\ln u dvol.
\ea
There holds 
\ba \label{shannon}
\frac{d}{dt} S_D(u)=\int_M \frac{|\nabla u|^2}{u} dvol,
\ea
whence $S_D(u)$ is increasing along the heat flow (\ref{heat}). Note that this leads to the following 
inequality 
\ba
\int_M u\ln u dvol \ge \int_M \bar u \ln \bar u dvol
\ea
for any nonegative function $u$ on $M$, where $\bar u$ denotes the average of $u$. Now the scalar curvature 
$R$ satisfies a nonlinear heat equation along the Ricci flow or the volume-normalized Ricci flow, hence it is 
analogous to the above $u$ of the heat flow. 
Obviously, Hamilton's entropy is simply the negative multiple of the differential Shannon entropy of 
$R$, where $R$ is viewed as the probability density function. (Note that $\int_M R dvol$ is a constant 
along the volume-normalized Ricci flow by Gauss-Bonnet theorem.) The  formula for Hamilton's entropy $S_0$ analogous to (\ref{shannon}) is  
\ba \label{H1}
\frac{dS_0}{dt}=-\int_M \frac{|\nabla R|^2}{R} dvol+\int_M R^2 dvol
\ea
along a smooth solution of the Ricci flow with positive scalar curvature. Indeed, since 
$\frac{\partial}{\partial t}R=\Delta R+R^2$ and $\frac{\partial}{\partial t}dvol=
-R dvol$ we have 
\ba
\frac{dS_0}{dt}&=&\int_M (\Delta R+R^2) \ln R dvol+\int_M (\Delta R+R^2)dvol
-\int_M R^2 \ln R dvol \nonumber \\ 
&=& -\int_M \frac{|\nabla R|^2}{R} dvol+\int_M R^2 dvol. \nonumber 
\ea
Let $f$ be Hamilton's potential function which is defined by $\Delta f=R-r$ with 
$\int_M f dvol=0$, where $r$ denotes the average of $R$. Then (\ref{H1}) implies
\ba \label{H3}
\frac{dS_0}{dt} =-\int_M \frac{|\nabla R+R \nabla f|^2}{R}dvol-2 \int_M |H|^2 dvol+8 \pi r,
\ea
where $H$ is the trace-free part of the Hessian of $f$. (This follows from 
(\ref{H1}), the formulas in the proof of Proposition 5.39 in [CK], and Gauss-Bonnet theorem.) 
The term $8 \pi r$ can be given explicitly
\ba
8\pi r=\frac{64 \pi^2}{V_0-8 \pi t}= -8 \pi \frac{d}{dt} \ln (V_0-8 \pi t), 
\ea
where $V_0$ is the volume of the initial metric. Hence we obtain the following monotonicity formula. 

\begin{theo} There holds along a smooth solution of the Ricci flow on $M$
\ba \label{H4}
\frac{d}{dt}[S_0+8 \pi \ln (V_0-8 \pi t)]=-\int_M \frac{|\nabla R+R \nabla f|^2}{R}dvol-2 \int_M |H|^2 dvol.
\ea
\end{theo}

In [H1], a simpler monotonicity formula is given for $S_0$ along a smooth solution of the 
volume-normalized Ricci flow. Indeed, there holds, similar to 
(\ref{H1}),  
\ba \label{H5}
\frac{dS_0}{dt}=-\int_M \frac{|\nabla R|^2}{R} dvol+\int_M (R-r)^2 dvol
\ea
along a smooth solution $g=g(t)$ of the volume-normalized Ricci flow with positive scalar curvature, 
which  leads to the monotonicity formula of Hamilton
\ba \label{H6}
\frac{dS_0}{dt} =-\int_M \frac{|\nabla R+R \nabla f|^2}{R}dvol-2 \int_M |H|^2 dvol,
\ea
where 
$f$ and $H$ are defined in the same way as the $f$ and $H$ in (\ref{H3}). This monotonicity 
formula immediately yields an upper bound for $S_0$, which is combined in [H1] with 
the differential Harnack inequality for $R$ to produce an upper bound for $R$. This is the 
key step in [H1] for establishing the convergence of the volume-normalized Ricci flow in the case 
that the initial metric has positive scalar curvature. (For another approach to the Ricci flow on the 
2-sphere based on the parabolic moving plane method in [Y1] we refer to [BSY]. 
\footnote{In [Y1], we 
neglected to mention Th.~Aubin's important contribution to the Yamabe problem.})

Next we consider 
a more general entropy quantity, the {\it adjusted entropy} $S_a$ 
\ba
S_a(g)=\int_M (R-a) \ln (R-a) dvol
\ea
for a given parameter $a$, assuming $R \ge a$. 
For a smooth solution  $g=g(t)$ of the volume-normalized Ricci flow we can choose 
$a=a(t)$ to be the solution of the ODE
\ba
\frac{da}{dt}=a(a-r)
\ea
and consider $S_{a(t)}(g(t))$. This ODE corresponds to the evolution equation for 
the scalar curvature 
\ba
\frac{\partial R}{\partial t}=\Delta R+R(R-r)
\ea
associated with the Ricci flow.
The quantity $S_{a(t)}(g(t))$ is precisely the modified entropy introduced by B.~Chow in [Ch].
(In the case of the Ricci flow, we can choose $a=a(t)$ to be the solution of the ODE $\frac{d}{dt}a=a^2$,
which corresponds to the evolution equation $\partial R/\partial t=\Delta R+R^2$.)
Although $S_{a(t)}(g(t))$ is not necessaily monotone, Chow is able to obtain an upper bound for it. 
He uses this upper bound together with other tools to show that the scalar curvature always becomes positive in finite time 
regardless what initial metric is given, see [Ch] or [CK] for details.   

Hamilton's entropy $S_0$ can obviously be defined on manifolds of dimensions $n \ge 3$, but no 
monotonicity formula for it has been found there. Note that part of the analogy of the differential 
Shannon entropy is lost there, because in general the integral of the scalar curvature  no longer stays a constant along 
the  volume-normalized Ricci  flow.
We think that a deeper understanding of Hamilton's entropy should be 
pursued.

\sect{Perelman's entropy functional} 

G.~Perelman introduced in [P1] the entropy functional ${\mathcal W}(g, f, \tau)$ and established its 
monotonicity along the Ricci flow. Let $M$ be a compact manifold of dimension $n \ge 2$. 
The definition of ${\mathcal{W}}(g, f, \tau)$ is as follows
\ba \label{entropy}
{\mathcal W}(g, f, \tau)=\int_M \left[ \tau(R+|\nabla f|^2)+f-n\right] \frac{e^{-f}}{(4\pi\tau)^{\frac{n}{2}}} dvol,
\ea
where $\tau$ is a positive number, $g$ is a Riemannian metric on $M$,  and $f\in C^{\infty}(M)$ satisfies 
\ba \label{vol-1}
\int_M \frac{e^{-f}}{(4\pi\tau)^{\frac{n}{2}}} dvol=1.
\ea
All geometric quantities in (\ref{entropy}) and (\ref{vol-1}) are associated with $g$. 

\begin{theo} (Perelman) Let $g=g(t)$ be a smooth solution of the Ricci flow on $M \times I$ 
for some interval $I$. Let $\tau=\tau(t)$ be a scalar function on $I$ with $\tau'=-1$, and 
$f=f(t)$ a smooth solution of the equation 
\ba
\label{nonlinearconjugate}
\frac{\partial f}{\partial t}=-\Delta f+|\nabla f|^2-R+\frac{n}{2\tau}
\ea
associated with $g=g(t)$ on $M \times I$. Then there holds
\ba
\frac{d {\mathcal W}}{d t}=2\tau \int_M |Ric+\nabla^2 f-\frac{1}{2\tau} g|^2 \frac{e^{-f}}{(4\pi\tau)^{\frac{n}{2}}} dvol
\ge 0
\ea
in $I$,  where ${\mathcal W}={\mathcal W}(g(t), f(t), \tau(t))$. 
\end{theo}

For a detailed proof of this theorem we refer to [KL]. In [P1], Perelman provided 
an interpretation of his entropy in terms of statistical thermodynamics, which we briefly explain here. 
Consider the thermodynamic system (in the sense of analogy) discribed by the pair $(g, f)$ with $g$ denoting a Riemannian metric on $M$ and $f$ a smooth 
function on $M$.  By well-known formulas in statistical mechanics the entropy of the 
system is given by 
\ba
S=\beta <E>+\ln Z,
\ea
where $\beta$ is the inverse of the temperature, $<E>$ the expectation value of energy, and 
$Z$ the partition function.  The formulas for $Z$ and $E$ are
\ba
Z=\int e^{-\beta E} d \mu(E),\,\,\, 
<E>=-\frac{\partial}{\partial \beta} \ln Z,
\ea
with $\mu$ denoting the ``density of states" measure. Perelman states that the partition function is 
given by 
\ba
Z=\int_M (-f+\frac{n}{2})dm,
\ea
where 
\ba
dm=\frac{e^{-f}}{(4\pi\tau)^{\frac{n}{2}}} dvol.
\ea
We do not attempt to explain this result here.  Now let the pair $(g, f)$ evolve 
by the Ricci flow coupled with (\ref{nonlinearconjugate}), with $\tau$ as given in Theorem \ref{Perelman} 
playing the role of time.  Then a standard computation yields
\ba
S=-{\mathcal{W}}(g, f, \tau).
\ea
We would also like to mention the following obervation in [KL]
\ba
{\mathcal{W}}(g, f, \tau)=\frac{d}{d\tau} \left(\tau \int_M (f-\frac{n}{2})dm\right).
\ea
 
Next we present a natural motivation of Perelman's entropy functional in terms of the logarithmic 
Sobolev inequality on the euclidean space ${\bf R}^n$.   The logarithmic Sobolev inequality of L.~Gross states 
\ba \label{Gross}
\int u^2 \ln u^2 d\mu \le 2 \int |\nabla u|^2 d\mu,
\ea 
assuming $u\in W^{1,2}_{loc}({\bf R}^n)$ and
 $\int u^2d\mu=1$, where
$d\mu=(2\pi)^{-\frac{n}{2}}e^{-\frac{|x|^2}{2}}dx.$ (We omit the integration domain which is 
the entire ${\bf R}^n$.)
This can be proved by employing the ordinary Sobolev inequality with the optimal constant 
and the product structure of the Euclidean spaces. Setting $u=(2 \pi)^{-\frac{n}{4}} e^{-f/2}$ 
we can convert (\ref{Gross}) into the following formulation due to Perelman [P1]
\ba \label{Perelman}
\int (\frac{1}{2}|\nabla f|^2+f-n)\frac{e^{-f}}{(2\pi)^{\frac{n}{2}}}dx \ge 0,
\ea
provided that $f \in W^{1,2}_{loc}({\bf R}^n)$, $\int e^{-f}(2\pi)^{-n/2}dx=1$ and $\int |\nabla f|^2 
e^{-f}dx<\infty$. 
Obviously, the left hand side of (\ref{Perelman}) is precisely 
${\mathcal{W}}(g_{euc}, f, \frac{1}{2})$, where $g_{euc}$ is the euclidean metric. To see how 
the general value of the parameter $\tau$ enters into the play, we can consider the 
log gradient version of the logarithmic Sobolev inequality 
\ba \label{loggradient}
\int u^2 \ln u^2 dx \le \frac{n}{2} \ln \left[ \frac{2}{\pi n e } \int |\nabla u|^2 dx \right]
\ea
for $u \in W^{1,2}({\bf R}^n)$  with $\int u^2 dx =1$. Since $\ln s \le \sigma s-\ln \sigma-1$ 
for all $s>0$ and $\sigma>0$  we deduce from (\ref{loggradient}) 
\ba
\int u^2 \ln u^2 dx\le \frac{n}{2} \ln \frac{2}{\pi n e }+ \frac{n\sigma}{2} \int |\nabla u|^2 dx
-\frac{n}{2} \ln \sigma -\frac{n}{2}.
\ea
To simplify the formula we replace $\sigma$ by  $\frac{2}{n} \sigma$ and obtain 
\ba
\int u^2 \ln u^2 dx\le \sigma \int |\nabla u|^2 dx 
-\frac{n}{2} \ln \sigma  -n-\frac{n}{2} \ln \pi .
\ea  
This form of the logarithmic Sobolev inequality has been used extensively in the theory of 
ultracontractivity of symmetric Markov processes.  We can absorb the terms $-\frac{n}{2} \ln \sigma 
-\frac{n}{2} \ln \pi$ if we replace $u$ by $(\pi \sigma)^{-n/4} u$. We obtain
\ba
\int u^2 \ln u^2 \frac{1}{(\pi \sigma)^{\frac{n}{2}}}dx \le \sigma \int |\nabla u|^2  \frac{1}{(\pi \sigma)^{\frac{n}{2}}} dx
-n,
\ea
where $u$ satisfies 
\ba
\int \frac{u^2}{(\pi \sigma)^{\frac{n}{2}}} dx=1.
\ea

Finally, we set $u=e^{-f/2}$ i.e. $f=-\ln u^2$ and $\sigma=4\tau$ to deduce, assuming $u>0$ 
\ba
-\int f \frac{e^{-f}}{(4\pi \tau)^{\frac{n}{2}}} dx\le 
\tau \int |\nabla u|^2  \frac{e^{-f}}{(4\pi \tau)^{\frac{n}{2}}}dx -n,
\ea
provided that $\int e^{-f}(4\pi \tau)^{-n/2}dx=1.$
We summarize the above compuations in a theorem. 

\begin{theo} Let $f \in W^{1,2}_{loc}({\bf R}^n)$ and $\tau>0$ such that 
\ba
\int \frac{e^{-f}}{(4\pi \tau)^{\frac{n}{2}}}dx=1
\ea
and 
\ba
\int |\nabla f|^2 e^{-f} dx<\infty.
\ea
 Then there holds
\ba \label{new}
\int [\tau |\nabla f|^2 +f -n] \frac{e^{-f}}{(4\pi \tau)^{\frac{n}{2}}}dx \ge 0.
\ea
\end{theo}

The  formula in (\ref{new}) obviously suggests the general formula of ${\mathcal{W}}(g, f, \tau)$ on 
an arbitary manifold. Of course, one still needs to come up with the idea of 
replacing $|\nabla f|^2 $ by $|\nabla f|^2+R$ in the general formula. The appearance of $R$ in the 
formula is quite natural in view of the classical total scalar curvature functional which 
plays an important role for the study of Einstein metrics. The fact that $R$ should be put together with 
$|\nabla f|^2$ is also natural in view of the common scaling weight of $R$ and $|\nabla f|^2$.  \\

\sect{The log entropy functional and the log Sobolev constant}

Motivated by the log gradient version of the logarithmic Sobolev inequality, the concept of 
log entropy functional was introduced in [Y6]. Let $M$ be a closed manifold of dimension $n\ge 2$.  
 We define the {\it log entropy} functional as follows

\ba \label{log1}
{\mathcal Y}_0(g, u)=-\int_M u^2 \ln u^2 dvol +\frac{n}{2} \ln  \left( \int_M (|\nabla u|^2+ \frac{R}{4}u^2)dvol
\right),
\ea
 where $g$ is a smooth metric on $M$ and $u \in W^{1,2}(M)$ satisfies 
\ba \label{eigenpositive}
\int_M (|\nabla u|^2+ \frac{R}{4} u^2) dvol >0.
\ea
Here, all geometric quantities are associated with $g$. 
More generally, we define the {\it log entropy} functional with {\it remainder} $a$ as follows 
\ba
{\mathcal Y}_{a}(g, u)=-\int_M u^2 \ln u^2 dvol +\frac{n}{2} \ln  \left( \int_M (|\nabla u|^2+ \frac{R}{4}u^2)dvol+a
\right).
\ea

Finally, we   define  the  {\it adjusted log entropy} with remainder $a$ as follows
\ba
{\mathcal Y}_{a}(g, u, t)=-\int_M u^2 \ln u^2 dvol +\frac{n}{2} \ln  \left( \int_M (|\nabla u|^2+ \frac{R}{4}u^2)dvol
+a 
\right) +4at.
\ea
Obviously, ${\mathcal Y}_a(g, u)={\mathcal Y}_a(g, u, 0)$. \\

Now we consider a smooth solution $g=g(t)$ of the Ricci flow on $M \times [\alpha, T)$ for some 
$\alpha<T$, where $\alpha$ is finite.  Let $u=u(t)$ be a smooth positive solution of the backward evolution equation 
\ba \label{uequation}
\frac{\partial u}{\partial t}=-\Delta u+\frac{|\nabla u|^2}{u}+\frac{R}{2}u
\ea
which is derived from the equation (\ref{nonlinearconjugate}) upon setting $u=e^{-f/2}$.
Let $\lambda_0$ denote the first eigenvalue of the operator $-\Delta+\frac{R}{4}$. 
The following monotonicity result was obtained in [Y6].

\begin{theo} \label{entropymonotone} Assume that $a >-\lambda_0(g(\alpha))$. Then 
${\mathcal Y}_a(t) \equiv {\mathcal Y}_{a}(g(t), u(t), t)$ is nondecreasing. Indeed, we have 
\ba \label{monotone1}
\frac{d}{dt}{\mathcal Y}_{a} \ge \frac{n}{4\omega} 
\int_M |Ric-2 \frac{\nabla^2 u}{u}+2 \frac{\nabla u \otimes \nabla u}{u^2}
-\frac{g}{16\omega}|^2 u^2dvol,
\ea
where 
\ba
\omega=\omega(t)=a+\int_M (|\nabla u|^2+\frac{R}{4}u^2)dvol|_{t},
\ea
which is positive. 
\end{theo}

This theorem was proved by combining Perelman's entropy monotonicity formula with 
a minimizing procedure. To see its consequence on the behavior of the logarithmic Sobolev inequality 
along the Ricci flow, we define for each $a>-\lambda_0(g)$ the {\it logarithmic Sobolev constant with the
$a$-adjusted scalar curvature 
potential}
\ba
C_{S, log, a}(M, g)&=&\inf\{-\int_M u^2 \ln u^2 dvol +\frac{n}{2} \ln  \left( \int_M (|\nabla u|^2+ (\frac{R}{4}+a)u^2)dvol 
\right):  \nonumber \\
&& u \in W^{1,2}(M), \int_M u^2 dvol=1 \}.
\ea 
In other words, $C_{S, log, a}(M,g)$ is the optimal constant, i.e. the maximal possible constant, such that the logarithmic 
Sobolev inequality 
\ba
\int_M u^2 \ln u^2 dvol \le \frac{n}{2} \ln  \left( \int_M (|\nabla u|^2+ \frac{R}{4}u^2)dvol+a\right) 
-C_{S, log, a}(M,g)
\ea
holds true for all $u \in W^{1,2}(M)$ with $\int_M u^2 dvol=1$.  
The following monotonicity results were established in [Y6] as  applications of Theorem \ref{entropymonotone}.

\begin{theo}  \label{logconstantmonotone} The adjusted logarithmic Sobolev inequality improves along the Ricci flow.
More precisely, $C_{S, log, a}(M, g(t))+4at$ is nondecreasing along an arbitary 
smooth solution $g(t)$ of the Ricci flow on $M$, provided that $a$ is greater than 
the negative mutilple of the $\lambda_0$ of the initial metric.  In particular, the logarithmic Sobolev 
inequality improves along the Ricci flow, i.e. 
$C_{S, log, 0}(M, g(t))$ is nondecreasing along the Ricci flow, provided that 
$\lambda_0>0$ at the start.    
\end{theo}

\begin{theo} \label{zerocase} Assume $\lambda_0=0$ at the start.   Then either $g=g(t)$ is a gradient soliton, in which case 
the logarithmic Sobolev constant $C_{S, log, a}(M, g(t))$ is independent of $t$ for any given $a>0$; 
or $C_{S, log, 0}(M, g(t))$ is nondecreasing on $[\epsilon, T)$ for each $\epsilon>0$. 
\end{theo}

\sect{The logarithmic Sobolev and  Sobolev inequalities along the Ricci flow}

The following results on the logarithmic Sobolev inequality along the Ricci flow 
were obtained in [Y3]. Consider a compact manifold $M$ of dimension $n \ge 3$.  Let $g=g(t)$ be a smooth solution of the Ricci flow 
on $M \times [0, T)$ for some (finite or infinite) $T>0$ with a given initial metric 
$g(0)=g_0$.\\

\begin{theo} For each $\sigma>0$ and each $t \in [0, T)$  there holds
\ba \label{sobolevA}
\int_M u^2 \ln u^2 dvol &\le& \sigma \int_M (|\nabla u|^2 +\frac{R}{4} u^2)dvol  -\frac{n}{2} \ln \sigma 
+A_1(t+\frac{\sigma}{4})+A_2
\ea
for all $u\in W^{1,2}(M)$ with  $\int_M u^2 dvol=1$, 
where
\ba
A_1&=&\frac{4}{\tilde C_S(M, g_0)^2 vol_{g_0}(M)^{\frac{2}{n}}}-\min R_{g_0}, \nonumber \\
A_2&=& n\ln \tilde C_S(M,g_0)+
\frac{n}{2}(\ln n-1), \nonumber
\ea
and 
all geometric quantities  are associated with the metric $g(t)$ (e.g. the volume form $dvol$ and the scalar curvature $R$), except the scalar curvature $R_{g_0}$, the modified Sobolev constant $\tilde C_S(M,g_0)$ (see Section 2 for its definition) and the volume $vol_{g_0}(M)$ which 
are those of the initial metric $g_0$. 
\end{theo}

\begin{theo} \label{positive} Assume that $\lambda_0(g_0) > 0$. For each $t\in [0, T)$ and each $\sigma>0$ there holds 
\ba \label{sobolevC}
\int_M u^2 \ln u^2 dvol \le \sigma \int_M (|\nabla u|^2 +\frac{R}{4}u^2) dvol
-\frac{n}{2}\ln \sigma +C
\ea
for all $u \in W^{1,2}(M)$ with $\int_M u^2 dvol=1$, where $C$ depends only on the dimension $n$, a positive lower bound for $vol_{g_0}(M)$, a nonpositive 
lower bound for $R_{g_0}$, an upper bound for $C_S(M, g_0)$, and a positive lower bound for 
$\lambda_0(g_0)$.
\end{theo}

The 2-dimensional case and the case $\lambda_0(g_0)=0$ are treated in [Y4] and 
[Y5]. The class of Riemannian manifolds $(M, g_0)$ with $\lambda_0(g_0) \ge 0$   is a very large one
and particularly significant from a geometric point of view. 
On the other hand, we  would like to point out that in general the assumption $\lambda_0(g_0) \ge 0$
indispensible in (\ref{positive}).  In other words, a uniform logarithmic 
Sobolev inequality like (\ref{sobolevC}) without the assumption $\lambda_0(g_0) \ge 0$ 
is false in general. Indeed, by [HI] there are  smooth solutions of the Ricci flow 
on torus bundles over the circle which exist for all time, have  bounded curvature, 
and collapse as $t\rightarrow \infty$. In view of the $\kappa$-noncollapsing estimate 
implied by (\ref{sobolevC}) (see Section 6)  a  uniform logarithmic Sobolev inequality like (\ref{sobolevC}) fails to hold along these solutions.

 Quantitative improvements of the above logarithmic Sobolev inequalities follow easily from the monotonicity 
theorems on the log Sobolev constant in the last section. Although the log gradient version of the 
logarithmic Sobolev inequality appears to be stronger than the above $\sigma$-version of 
the logarithmic Sobolev inequality, they are actually equivalent. On the other hand, 
the $\sigma$-version is often more convenient for applications.  
In particular, the ordinary Sobolev inequalities can be derived from it, which we address next.

It has been known for some time that there exist close relations between the logarithmic 
Sobolev inequality, the $W^{1,2}$ Sobolev inequality and the so-called ultracontractivity 
of the heat semigroup of the associated Schr\"{o}dinger operator.  Indeed, they are equivalent.
This theory is presented e.g. in [D] in the general and abstract set-up of 
symmetric Markov processes.  Besides PDE arguments, basic spectral analysis of self-adjoint operators and 
some basic theorems of harmonic analysis 
such as the Riesz-Thorin interpolation theorem and the Marcinkiewicz  interpolation theorem play a crucial 
role. For the purpose of precise geometric estimates, transparent presentation  and additional implications, the theory in [D] 
is adapted in [Y3] to the geometric set-up there and worked out in complete and self-contained details. 
(The paper [Z] provided help for us to find the reference [D].)
Based on this theory,  in [Y3] the following results on the Sobolev inequalties along the Ricci flow are derived from the 
above logarithmic Sobolev inequalities. 
Let $g=g(t)$ be a smooth solution of the Ricci flow on $M \times [0, T)$ with a given initial metric $g_0$ as 
before. 

\begin{theo} \label{sobolevDD} Assume that $\lambda_0(g_0) > 0$.  There is a  positive constant $A$ depending only on 
the dimension $n$, a nonpositive lower bound for $R_{g_0}$, a positive lower bound for $vol_{g_0}(M)$, an upper bound for $C_S(M,g_0)$,
and a positive lower bound for $\lambda_0(g_0)$, such that for each $t \in [0, T)$ and 
all $u \in W^{1,2}(M)$ there holds 
\ba \label{sobolevD}
\left( \int_M |u|^{\frac{2n}{n-2}} dvol \right)^{\frac{n-2}{n}} \le 
A\int_M (|\nabla u|^2+\frac{R}{4}u^2) dvol,
\ea
where all geometric quantities except $A$  are associated with $g(t)$. 
\end{theo}

\begin{theo} \label{sobolevDDD} Assume $T<\infty$. There are positive constants $A$ and $B$ depending only on the dimension $n$, 
a nonpositive lower bound for $R_{g_0}$, a positive lower bound for $vol_{g_0}(M)$, an upper bound for $C_S(M,g_0)$,
and an upper bound for $T$, such that for each $t \in [0, T)$ and 
all $u \in W^{1,2}(M)$ there holds 
\ba \label{sobolevD*}
\left( \int_M |u|^{\frac{2n}{n-2}} dvol \right)^{\frac{n-2}{n}} \le 
A\int_M (|\nabla u|^2+\frac{R}{4}u^2)dvol+B \int_M u^2 dvol,
\ea
where all geometric quantities except $A$ and $B$  are associated with $g(t)$. 
\end{theo}

\sect{The $\kappa$-noncollapsing estimates}

We first recall the following definitions [P1] [Y2].\\

\noindent {\bf Definition 1} Let $g$ be a Riemannian metric on a manifold $M$ of dimension
$n$. Let $\kappa$ and $\rho$ be positive numbers. We say that $g$ is {\it $\kappa$-noncollapsed on the scale} 
$\rho$,
if $g$ satisfies $vol(B(x, r))\ge \kappa r^n $
for all $x \in M$ and $r>0$ with the properties $r < \rho$ and $\sup\{|Rm|(x): x\in B(x, r)\}
\le r^{-2}$.  We say that a family of Riemannian metrics $g=g(t)$ is $\kappa$-noncollapsed on the scale 
$\rho$, if $g(t)$ is $\kappa$-noncollapsed on the scale $\rho$ for each $t$ (in the given domain). 
\\

\noindent {\bf Definition 2} Let $g$ be a Riemannian metric on a manifold $M$ of dimension
$n$. Let $\kappa$ and $\rho$ be positive numbers. We say that $g$ is {\it $\kappa$-noncollapsed on the scale 
$\rho$ relative to upper bounds of the 
scalar curvature},
if $g$ satisfies  $vol(B(x, r))\ge \kappa r^n $
for all $x \in M$ and $r>0$ satisfying $r < \rho$ and $\sup\{R(x): x\in B(x, r)\}
\le r^{-2}$. 
\\

The $\kappa$-noncollapsing property of the Ricci flow on a compact manifold can be 
derived directly from Perelman's entropy monotonicity, as done in [P1] and 
[Y2], see also [KL].  But the Sobolev inequalitities along the Ricci flow obtained in 
[Y3][Y4] and [Y5] and as presented in the last section lead to better results in several ways.
First, we obtain explicit $\kappa$-noncollapsing estimates which have clear and 
rudimentary geometric dependences on the initial metric. Second, the 
estimates are uniform up to $t=0$, which is not immediately clear from the direct 
arguements. Third,  the estimates hold true uniformly for all time (up to infinity), 
provided that $\lambda_0(g_0) \ge 0$.  We present in this section only the results from [Y3]. 
First we have the following general result on the $\kappa$-noncollapsing 
estimate implied by the Sooblev inequality.

\begin{theo} \label{sobolevDDDD} Consider the Riemannian manifold $(M,g)$ for a given metric $g$, such that 
for some $A>0$ and $B>0$ the Sobolev inequality 
\ba \label{lemmasobolev}
\left( \int_M |u|^{\frac{2n}{n-2}} dvol \right)^{\frac{n-2}{n}} \le 
A\int_M (|\nabla u|^2+\frac{R}{4}u^2) dvol +B\int_M u^2 dvol
\ea
holds true for all $u \in W^{1,2}(M)$. 
Let $L>0$. Assume $R\le \frac{1}{r^2}$ on a geodesic ball $B(x, r)$ with $0<r\le L$. Then 
there holds 
\ba \label{noncollapse1}
vol(B(x, r)) \ge \left(\frac{1}{2^{n+3}A+2BL^2}\right)^{\frac{n}{2}} r^n.
\ea
\end{theo} 

Basically, the Sobolev inequality implies a Faber-Krahn inequality for 
the first Dirichlet eigenvalue of subdomains in $B(x, r)$ under the assumption 
$R \le \frac{1}{r^2}$ on $B(x, r)$. By an elementary iteration argument in [Ca] 
we then arrive at the desired volume estimate.    

Combining Theorem \ref{sobolevDD}  and Theorem \ref{sobolevDDD} with Theorem \ref{sobolevDDDD} 
we arrive at the following $\kappa$-noncollapsing estimates in [Y3]. 

\begin{theo} \label{noncollapseI} Assume that $\lambda_0(g_0) > 0$.  Let $t \in [0, T)$. Consider the Riemannian manifold 
$(M, g)$ with $g=g(t)$. Assume $R\le \frac{1}{r^2}$ on a geodesic ball $B(x, r)$ with $r>0$. Then 
there holds 
\ba \label{noncollapse}
vol(B(x, r)) \ge \left(\frac{1}{2^{n+3}A}\right)^{\frac{n}{2}} r^n,
\ea
where $A$ is from Theorem \ref{sobolevDD}. In other words, the flow $g=g(t), t\in [0, T)$ is $\kappa$-noncollapsed relative to 
upper bounds of the scalar curvature on all scales. 
\end{theo}

\begin{theo} \label{noncollapseII} Assume that $T<\infty$.  Let $L>0$ and $t \in [0, T)$. Consider the Riemannian manifold 
$(M, g)$ with $g=g(t)$. Assume $R\le \frac{1}{r^2}$ on a geodesic ball $B(x, r)$ with $0<r \le L$. Then 
there holds 
\ba \label{noncollapse*}
vol(B(x, r)) \ge \left(\frac{1}{2^{n+3}A+2BL^2}\right)^{\frac{n}{2}} r^n,
\ea
where $A$ and $B$ are from Theorem \ref{sobolevDDD}. 
\end{theo} 

The imorptance of  $\kappa$-noncollapsing estimates is that they yield estimates for the 
injectivity radius under the assumption of bounds for curvatures.  In particular, they 
enable one to obtain smooth blow-up limits for the Ricci flow, which is crucial for blow-up analysis of 
singularities of the Ricci flow.   The said estimates for the injectivity radius is implied 
by the result in [CGT] on estimates of the injectivity radius. We state here the local formulation of this 
result as given in [Y1]. 

\begin{theo}     Let $(M, g)$ be a Riemannian manifold of dimension $n$.
Assume that the sectional curvatures $K_g$ of $g$ satisfies 
$\kappa_1 \le K_g \le \kappa_2$ on a geodesic ball $B(p, r_0)$ in $(M, g)$, 
such that $r_0 \le d(p, \partial M)$. (For the general definition of 
$d(p, \partial M)$ see [Y2]. Note that $d(p, \partial M)=\infty$ if 
$M$ is closed. )  Set $r_1=\frac{1}{4}\min\{r_0, \frac{\pi}{4\sqrt{\kappa_2}}\}$. 
Then the injectivity radius $i(q)$ at any $q 
\in B(p, r_1)$ 
satisfies 
\be 
i(q) \ge r_2, \tag{B.1} 
\ee
where 
\be
r_2=\frac{r_1}{2} \left(1+\frac{V_{\kappa_1}(2r_1)^2}{vol_g(B(p, r_1))V_{\kappa_1}(r_1)}\right)^{-1}, \tag{B.2}
\ee
 $r_1=\frac{1}{4}\min\{r_0, \frac{\pi}{\sqrt{\kappa_2}}\}$,  and for any $r>0$, $V_{\kappa_1}(r)$ denotes the 
volume of a geodesic ball of radius $r$ in the $n$-dimensional model space (a simply connected 
complete Riemannian manifold) of sectional curvature $\kappa_1$. 
\end{theo}

\sect{The modified Ricci flows} 

By scaling invariance, the above results extend to  
the modified Ricci flow 
\ba \label{modified}
\frac{\partial g}{\partial t}=-2Ric+\lambda(g, t) g
\ea
with a smooth scalar function $\lambda(g, t)$ independent of $x \in M$. 
The volume-normalized 
Ricci flow 
\ba \label{volumenormalize}
\frac{\partial g}{\partial t}=-2Ric+\frac{2}{n} {\hat R} g
\ea
on a closed manifold, with $\hat R$ denoting the average scalar curvature, is  an example of the modified 
Ricci flow. The $\lambda$-normalized Ricci flow 
\ba \label{lambda}
\frac{\partial g}{\partial t}=-2Ric+\lambda g
\ea
for a constant $\lambda$ is another example. (Of course, it reduces to the Ricci flow when 
$\lambda=1$.) The normalized K\"{a}hler-Ricci flow is a special 
case of it. 

We have e.g.~the following results from [Y3].
\\

\begin{theo}  \label{modifyI} Theorem \ref{sobolevDD} and Theorem \ref{noncollapseI} extend to 
the modified Ricci flow.   
\end{theo}

Let $g=g(t)$ be a smooth solution of the modified Ricci flow 
(\ref{modified}) on $M \times [0, T)$ for some (finite or infinite) $T>0$,  with 
a given initial metric $g_0$. We set 
\ba
T^*=\int_0^T e^{-\int_0^t \lambda(g(s), s)ds} dt.
\ea

\begin{theo} \label{modifyII} Assume that $T^*<\infty$.  \\
1) There are positive constants $A$ and $B$ depending only on the dimension $n$, 
a nonpositive lower bound for $R_{g_0}$, a positive lower bound for $vol_{g_0}(M)$, an upper bound for $C_S(M,g_0)$,
and an upper bound for $T^*$, such that for each $t \in [0, T)$ and 
all $u \in W^{1,2}(M)$ there holds 
\ba \label{sobolevH}
\left( \int_M |u|^{\frac{2n}{n-2}} dvol \right)^{\frac{n-2}{n}} \le 
A\int_M (|\nabla u|^2+\frac{R}{4}u^2)dvol+Be^{-\int_0^t \lambda(g(s), s)ds} \int_M u^2 dvol. 
\ea
2)  Let $L>0$ and $t \in [0, T)$. Consider the Riemannian manifold 
$(M, g)$ with $g=g(t)$. Assume $R\le \frac{1}{r^2}$ on a geodesic ball $B(x, r)$ with $0<r \le L$. Then 
there holds 
\ba \label{noncollapseH}
vol(B(x, r)) \ge \left(\frac{1}{2^{n+3}A+2Be^{-\int_0^t \lambda(g(s), s)ds}L^2}\right)^{\frac{n}{2}} r^n.
\ea
\end{theo}

Combining Theorems \ref{modifyI} and \ref{modifyII}  with Perelman's scalar curvature estimate [ST] we obtain the 
following result from [Y3]. \\

\begin{theo} Let $g=g(t)$ be a smooth solution of the 
normalized K\"{a}hler-Ricci flow 
\ba
\frac{\partial g}{\partial t}=-2Ric+2\gamma g
\ea
on $M \times [0, \infty)$ with a positive first Chern class, where $\gamma$ is the positive 
constant such that the Ricci class equals $\gamma$ times the K\"{a}hler class.   (We assume that 
$M$ carries such a K\"{a}hler structure.) Then 
the Sobolev inequality (\ref{sobolevH})  holds  true with $\lambda(g(s), s)=2\lambda$. 
Moreover, there is a positive constant $L$ depending only on the initial metric $g_0=g(0)$ and the dimension $n$ such that 
the inequality (\ref{noncollapseH}) holds true for all $t \in [0, T)$ and $0<r \le L$.

If $\lambda_0(g_0)>0$, then  
the Sobolev inequality (\ref{sobolevD}) holds true for $g$.  Moreover, there is a  positive constant 
depending only on the initial metric $g_0$ and the dimension $n$ such that 
the inequality (\ref{noncollapse}) holds true for all $t \in [0, T)$ and $0<r \le L$. 
Consequently, blow-up limits of $g$ at the time infinity  
satisfy (\ref{noncollapse}) for all $r>0$ and 
the Sobolev inequality 
\ba
\label{sobolevH}
\left( \int_M |u|^{\frac{2n}{n-2}} dvol \right)^{\frac{n-2}{n}} \le 
A\int_M |\nabla u|^2 dvol
\ea
for all $u$. 
(In particular, they must be noncompact.)  
\end{theo}

\sect{Furthe Sobolev inequalities}

As is well-known, the case $p=2$ of the $W^{1,p}(M)$ Sobolev inequalities is used most often in 
analysis and geometry. However, it is of high 
interest to understand the situation $1<p<2$ and $2<p<n$, both from the point of view of 
a deeper understanding of the theory and the point of view of further applications.  
In [Y7], $W^{1,p}$ and $W^{2,p}$ Sobolev inequalities for general $p$ along the Ricci flow 
were derived using several different methods, including Bessel potentials and 
Riesz transforms.  Indeed,  general results on deriving 
further Sobolev inequalities from a given Sobolev inequality were obtained in [Y7].   
We state here part of the results from [Y7].

\begin{theo} \label{highI} Assume that $\lambda_0(g_0) > 0$. Let $2<p<n$.  
Then there holds for each $t \in [0, T)$ and 
all $u \in W^{1,p}(M)$  
\ba \label{sobolevD}
\left( \int_M |u|^{\frac{np}{n-p}} dvol \right)^{\frac{n-p}{n}} \le 
A\left[(\max R^++1)vol(M)^{\frac{2}{n}}\right]^{\frac{m(p) p}{2}} \int_M (|\nabla u|^p+|u|^p) dvol,
\nonumber \\
\ea
where all geometric quantities are associated with $g(t)$,  except the constant $A$, 
which can be bounded from above in terms of the dimension $n$, a nonpositive lower bound for $R_{g_0}$, a positive lower bound for $vol_{g_0}(M)$, an upper bound for $C_S(M,g_0)$,
a positive lower bound for $\lambda_0(g_0)$, and an upper bound for $\frac{1}{n-p}$. The quantity 
$(\max R^++1)vol(M)^{\frac{2}{n}}$ is at time $t$ and the 
number $m(p)$ is defined as follows: $m(p)=2^{k+1}$ for 
$p \in (p_k, p_{k+1}]$, $p_0=2$, and $p_{k+1}=\frac{n^2 p_k}{(n-p_k)^2+np_k} $ for 
$k \ge 0$. 
\end{theo}

\begin{theo} \label{highII}  Assume $T<\infty$. 
Let $2<p<n$.  
Then there holds for each $t \in [0, T)$ and 
all $u \in W^{1,2}(M)$  
\ba \label{sobolevD*}
\left( \int_M |u|^{\frac{np}{n-p}} dvol \right)^{\frac{n-p}{n}} \le 
A\left[ 1+(\max R^++1)vol(M)^{\frac{2}{n}}\right]^{\frac{m(p)p}{2}} \int_M (|\nabla u|^p+|u|^p)dvol,
\nonumber \\
\ea
where all geometric quantities are associated with $g(t)$,  except the constant 
$A$, which  can be bounded from above in terms of the dimension $n$, 
a nonpositive lower bound for $R_{g_0}$, a positive lower bound for $vol_{g_0}(M)$, an upper bound for $C_S(M,g_0)$,
an upper bound for $T$, and an upper bound for $\frac{1}{n-p}$. 
The quantity 
$(\max R^++1)vol(M)^{\frac{2}{n}}$ is at time $t$ and the 
number $m(p)$ is the same as in Theorem \ref{highI}. 
\end{theo}

These two theorems are proved by employing an induction scheme based on the 
H\"{o}lder inequality.  Note that this scheme cannot be extended to $1<p<2$. 
Next we have the results on nonlocal Sobolev inequalities in terms of the 
$(1,p)$-Bessel norm, which is defined as follows
\ba
\|u\|_{B, 1, p}= \|(-\Delta+1)^{\frac{1}{2}}u\|_p
\ea
for $u \in W^{1,p}(M)$.

\begin{theo} Assume that $\lambda_0(g_0)>0$.  Let $1<p<n$.   There is a  positive constant $C$ depending only on 
the dimension $n$, a positive lower bound for $\lambda_0(g_0)$, a positive lower bound for $vol_{g_0}(M)$, an upper bound for $C_S(M,g_0)$, an upper bound for 
$\frac{1}{p-1}$, and an upper bound for $\frac{1}{n-p}$, such that for each $t \in [0, T)$ and 
all $u \in W^{1,p}(M)$ there holds 
\ba  \label{SB2-p}
\|u\|_{\frac{np}{n-p}} \le C (1+R_{max}^+)^{\frac{1}{2}}\|u\|_{B,1,p}.
\ea
\end{theo}

\begin{theo} Assume $T<\infty$ and $1<p<n$.  There is a positive constant $C$ depending only on the dimension $n$, 
a nonpositive lower bound for $R_{g_0}$, a positive lower bound for $vol_{g_0}(M)$, an upper bound for $C_S(M,g_0)$,
an upper bound for $T$, an upper bound for $\frac{1}{p-1}$, and an upper bound 
for $\frac{1}{n-p}$,  such that for each $t \in [0, T)$ and 
all $u \in W^{1,p}(M)$ there holds 
\ba \label{SB3-p}
\|u\|_{\frac{np}{n-p}} \le C (1+ R_{max}^+)^{\frac{1}{2}} \|u\|_{B,1,p}.
\ea
\end{theo} 
 
To convert the above results into ordinary Sobolev inequalities, the method of  Riesz transforms 
was used in [Y7]. Based on the $L^p$ estimates for the Riesz transform due to D.~Bakry [B] the following 
results were obtained in [Y7].

\begin{theo} Assume $\lambda_0(g_0)>0$.  Let $1<p<n$.   There is a  positive constant $C$ depending only on 
the dimension $n$, a positive lower bound for $\lambda_0(g_0)$, a positive lower bound for $vol_{g_0}(M)$, an upper bound for $C_S(M,g_0)$, an upper bound for 
$\frac{1}{p-1}$, and an upper bound for $\frac{1}{n-p}$, such that for each $t \in [0, T)$ and 
all $u \in W^{1,p}(M)$ there holds 
\ba 
\|u\|_{\frac{np}{n-p}} \le C (1+R_{max}^+)^{\frac{1}{2}} (\|\nabla u\|_{p}+(1+\kappa)\|u\|_p).
\ea
\end{theo}

\begin{theo} Assume $T<\infty$ and $1<p<n$.  There is a positive constant $C$ depending only on the dimension $n$, 
a nonpositive lower bound for $R_{g_0}$, a positive lower bound for $vol_{g_0}(M)$, an upper bound for $C_S(M,g_0)$,
an upper bound for $T$, an upper bound for $\frac{1}{p-1}$, and an upper bound 
for $\frac{1}{n-p}$,  such that for each $t \in [0, T)$ and 
all $u \in W^{1,p}(M)$ there holds 
\ba \label{Gallot}
\|u\|_{\frac{np}{n-p}} \le C (1+R_{max}^+)^{\frac{1}{2}} (\|\nabla u\|_{p}+(1+\kappa)\|u\|_p).
\ea
\end{theo} 

Similar results involving an integral norm of the Ricci curvature were also obtained in [Y7].
One should compare the above results with Gallot's estimate of the isoperimetric constant [G]
which implies an estimate for the Sobolev inequalities. In contrast to Gallot's estimates, no upper bound
for the diameter nor positive lower bound for the volume of $g(t)$ is assumed. Moreover, 
in (\ref{Gallot}) the lower bound for the Ricci curvature does not appear in front 
of $\|\nabla u\|_p$.  (Gallot's estimates lead to Sobolev inequalities in which  
the lower bound for the Ricci curvature gets involved with $\|\nabla u\|_p$.)

\sect{The Ricci flow with surgeries} 

The key construction in Perelman's work on the Ricci flow and geometrization of 3-manifolds is 
the Ricci flow with surgeries [P2], which extends Hamilton's earlier work on surgeries of the Ricci flow 
in a substantial way.  (For very nice accounts of Perelman's work we refer to 
[MT] and [KL].) Perelman's work provides a very clear picture of the behavior of 
the Ricci flow near blow-up singularities. Indeed, regions of large curvature are classified 
topologically and geometrically,  and 
one can find nice portions of the manifold near blow-up singularities where surgeries can be 
performed.  

Let $g=g(t)$ be a smooth solution of the Ricci flow on a closed 3-manifold $M$. 
(For technical reasons one needs to normalize the initial metric in a natural way.)
Assume that $g(t)$ is smooth on $[0, T)$ for a finite $T$, but becomes singular somewhere as 
$t \rightarrow \infty$.  (The situation of a Ricci flow with surgeries becoming singular 
upon approaching a finite time is similar.)  By Perelman's work, we can find a maximal region 
$\Omega \subset M$ such that $g(t)$ converges smoothly in $\Omega$ to a metric 
$g_T$.  Then the curvature of $g_T$ blows up as one approaches the boundary of 
$\Omega$. For  $\rho>0$ we consider the region $\Omega_{\rho}=\{ 
x \in \Omega: R < \rho^{-2}\}$ defined in terms of $g_T$. For sufficiently small $\rho$, $\partial \Omega_{\rho}$ consists 
of boundaries of $\epsilon$-necks whose interiors lie outside of $\Omega_{\rho}$. 
An $\epsilon$-neck is, after a suitable rescaling, $\epsilon$-close to 
the product $S^2 \times [-\epsilon^{-1}, \epsilon^{-1}]$ in the 
$C^{1/\epsilon}$ topology.  Let $Z_0$ be one of such $\epsilon$-necks and 
$Z$ the component of $\Omega-\Omega_{\rho}$ containing $Z_0$. If $Z$ is compact, 
then it is an $\epsilon$-cap, an $\epsilon$-neck or $\epsilon$-tube. In the last two cases, the  both boundary components of $Z$ are
contained in $\partial \Omega_{\rho}$.  If $Z$ is noncompact, then it is an $\epsilon$-horn 
consisting of   $\epsilon$-necks of increasing magnitudes of curvature, with the curvature going to 
infinity at the end.  We'll call $Z$ an {\it attached $\epsilon$-horn} for $\Omega_{\rho}$. 

 The surgery will be performed on the attached
$\epsilon$-horns for $\Omega_{\rho}$.  Indeed, for each attached $\epsilon$-horn $Z$, we pick a suitable 
$\epsilon$-neck in $Z$ and carry out the surgery on it.  The surgery consists of 
cutting the $\epsilon$-neck and hence $Z$ into two parts, throwing out the noncompact part of $Z$, and gluing in a sufficiently long capped 
cylinder.  The original metric $g_T$ will be interpolated inside of the $\epsilon$-neck with a standard metric on the capped 
cylinder.  After the surgery we restart the Ricci flow with the  
metric resulting from the surgery as the initial metric. This yields the first stage of 
a Ricci flow with surgeries.  Repeating the process whenever we run into singularities at 
a finite time
we then 
obtain a Ricci flow with surgeries on its maximal time interval.  

The surgeries have to be done in a way such that the following three key properties of the Ricci flow 
are preserved for the extended solution of the Ricci flow after the surgery: 1) the Hamilton-Ivey pinching, 2) the $\kappa$-noncollapsing property, and 3) 
the canonical neighborhood property.  The first property is the easiest to preserve. Indeed, 
a trick of Hamilton in [H2] can be applied without any difficulty.  In [P2], Perelman uses the reduced length and the reduced 
volume to handle the second property. Because of the surgery, the behavior of the ${\mathcal L}$-geodesics 
are rather complicated, and the arguments are also forced to be rather complicated. Furthermore, the second and third properties are established in [P2] by
rather involved combined arguments.

In [Y8], the arguments in [P2] for preserving the $\kappa$-noncollapsing property are replaced by  
a much simpler and more transparent argument . Indeed, the Sobolev inequality is first established for the 
metric resulting from the surgery. Then the Sobolev inequalty will continue to hold 
by the results in [Y3], as presented in Section 5. Finally, we obtain a $\kappa$-noncollapsing estimate 
as a consequence of the Sobolev inequality, as in [Y3] and explained in Section 6.
As an important feature of these arguments, the $\kappa$-noncollapsing property is 
established independent of the canonical neighborhood property, making the choice of the surgery 
parameters much simpler. It also becomes easier to establish the canonical neighborhood 
property.

The main results in [Y8] are as follows.
 
\begin{theo}
Let $n=3$ and $g=g(t)$ be a Ricci flow with surgeries as constructed 
in [P2] on its maximal time interval $[0, T_{max})$, with suitably chosen surgery parameters.  Let $g_0=g(0)$. Let $m(t)$ denote the number of 
surgeries which are performed up to the time $t\in (0, T_{max})$. Then there holds at  each $t \in [0, T_{max})$
\ba \label{sobolevJ*}
\left( \int_M |u|^{6} dvol \right)^{\frac{1}{3}} \le 
A(t)\int_M (|\nabla u|^2+\frac{R}{4}u^2)dvol+B(t) \int_M u^2 dvol
\ea
for all $u \in W^{1,2}(M)$, where $A(t)$ and $B(t)$ are bounded from above in terms of a nonpositive lower bound for $R_{g_0}$, a positive lower bound for $vol_{g_0}(M)$, an upper bound for $C_S(M,g_0)$, and
an upper bound for $t$.  

If $\lambda_0(g_0)>0$, then there holds at each $t \in [0, T_{max})$
\ba \label{sobolevJ}
\left( \int_M |u|^{6} dvol \right)^{\frac{1}{3}} \le 
A(t)\int_M (|\nabla u|^2+\frac{R}{4}u^2)dvol
\ea
for all $u \in W^{1,2}(M)$, where $A(t)$ is bounded from above in terms of  a nonpositive lower bound for $R_{g_0}$, a positive lower bound for $vol_{g_0}(M)$, an upper bound for $C_S(M,g_0)$,
a positive lower bound for $\lambda_0(g_0)$, and an upper bound for $m(t)$. 

$\kappa$-noncollapsing estimates follow from the above Sobolev inequalties.
\end{theo}

This theorem is based on a general result on  the Sobolev inequality on manifolds with surgeries and 
the specific patterns of surgeries for the Ricci flow in dimension 3 as discussed above.

\end{document}